\documentclass{amsart}
\usepackage{tabularx}
\usepackage{mathrsfs}
\usepackage{amsmath,amssymb,amsfonts,amsthm}
\usepackage{caption}
\usepackage{ctable}
\usepackage{geometry}
\usepackage{rotating}
\usepackage{makecell}
\usepackage{blkarray}

\title{Central Configurations with Dihedral symmetry}
\author{Tingjie Zhou}
\author{Zhihong Xia}
\address{Department of Mathematics, Southern University of Science and
 Technology, Shenzhen, China}
\address{Department of Mathematics, Northwestern University, Evanston,
  IL 60208 USA}
\email{11930530@mail.sustech.edu.cn, xia@math.northwestern.edu}
\date{\today}

\begin{document}

\maketitle

\begin{abstract}
  As an application of the representation theory for the dihedral
  groups, we study the symmetric central configurations in the
  $n$-body problem where $n$ equal masses are placed at the vertices
  of a regular $n$-gon. Since the Hessian matrices at these configurations
  are typically very large, particularly when $n$ is large, computations of
  their eigenvalues present a challenging problem. However, by
  decomposing the action of the dihedral groups into irreducible
  representations, we show that the Hessians can be simplified to a
  block-diagonal matrix with small blocks, of the sizes
  at most $2\times 2$. This is due to fact that the action of a
  dihedral group can be represented as a block-diagonal matrix with
  small irreducible blocks. In the end, the eigenvalues can be
  explicitly obtained by simply computing eigenvalues of these small block matrices.

\end{abstract}

\section{introduction}

For the Newtonian $n$ body problem in the plane, let $q_i$ and $m_i$ be the position and the mass for the $i$th particle for $i =1, 2, \ldots, n$. The central configurations are special configurations of the $n$ particles such that the relative distance of the particles can be maintained throughout a particular solution. It turns out that the central configurations are critical points of the function $\sqrt{I}U$, where
$$I=\frac{1}{2}\sum_{i=1}^{n}m_i|q_i|^2,\quad U=\sum_{1 \le i<j \le n} \frac{m_im_j}{|q_i-q_j|}.$$
The function $I$ is the moment of inertia of the $n$-body system and $U$
is the potential function.
Let $p_i=m_i\dot{q}_i$, the $n$-body system is Hamiltonian with Hamiltonian function
\begin{equation*}
    H=\sum_{i=1}^{n} \frac{|p_i|^2}{2m_i}-U.
\end{equation*}
The equation of motions are 
\begin{equation} \label{h1}
\begin{aligned}
\dot{q}_i=\frac{\partial H}{\partial p_i}=\frac{p_i}{m_i}, \\
\dot{p}_i=-\frac{\partial H}{\partial q_i}=\frac{\partial
  U}{\partial q_i}.
\end{aligned}
\end{equation}
A notable central configuration is where the $n$ equal masses are positioned at the vertices of a regular $n$-gon. This configuration exhibits inherent symmetries that offer avenues for simplification and deeper understanding. As it is a critical point for the function $\sqrt{I}U$, it is crucial to investigate the eigenvalues of the Hessian to uncover additional properties of the solution.

For the Newtonian $n$-body problem, for each planar central configuration $z_0 \in \mathbb{R}^{2n}$, there is correspondingly a class of periodic solutions where $n$ particles move in concentric circles, keeping their relative positions fixed. These solutions are called {\em relative equilibria}\/. A relative equilibria is degenerate if the Hessian for the corresponding central configuration has a non-trivial nullspace. In 1976, Palmore\cite{MR420713} studied the degeneracies of relative equilibria for the $n+1$-body problem, calculating eigenvalues by decomposing the tangent space into invariant subspaces. Through this approach, Palmore identified a critical mass value $m^*$ at which the relative equilibria degenerate. Building upon this work, Moeckel\cite{MR1350320} explored the linear stability of such symmetric configurations in the $n$-body problem. By identifying invariant subspaces and employing factorizations for characteristic polynomials, Moeckel further streamlined the analysis, emphasizing the crucial role of symmetry in this context.

 Leandro\cite{2017arXiv170502701L} systematically discussed the significance of symmetry in 2017, particularly focusing on ring systems. Leandro applied group representation theory and dihedral groups to exploit inherent symmetries, leading to a factorization of stability polynomials. Expanding on this methodology, Leandro\cite{MR3951830} applied these principles to investigate the structure and stability of the rhombus family of relative equilibria. These studies all revolve around computing eigenvalues for the Hessian, showcasing how symmetry aids in simplifying calculations.

In 2008, Xia\cite{doi:10.1063/1.2993622} introduced a novel method for studying symmetries within the Newtonian $n$-body problem. This approach was exemplified by determining eigenvalues for the Hessian of specific functions, showcasing how group representation and symmetry groups can simplify intricate problems. Expanding on Xia's work, we further refined this method in \cite{MR4355921}, applying it to streamline calculations for certain eigenvalue problems and elucidating the degenerate point $m^*$ in the $4$-body problem studied by Palmore\cite{MR420713}.

In this paper, we build on these foundations by presenting a new formulation of this method. By finding the irreducible representation for the dihedral groups, we can easily compute eigenspaces for these groups acting on $R^{2n}$. From these computations, we identify the smallest subspaces that's simultaneously invariant under the actions of all elements of the symmetry group. It turns out that the dimension of each of these subspaces is no larger than 2. We then decompose the whole space as a direct sum of these invariant subspaces. In the end, matrices that are invariant under the dihedral group action can be simplified into block diagonal forms with the size of each block at most $2\times 2$. We further apply this method to analyze the Newtonian $n$-body problem and derive explicit formulas for eigenvalues of the corresponding Hessians. In Section 2 we establish the framework, while in Section 3 we will apply the framework to the regular $n$-gon central configuration in several situations.

For a comprehensive overview of group representation and central configurations in the Newtonian $n$-body problem, we refer interested readers to \cite{MR4355921}.

\

The method we presented in this paper is applicable to general cases where a linear function is invariant under a group action.

\section{Block diagonal form}
\subsection{The irreducible representation of the dihedral groups} The dihedral group $D_n$ is the symmetry group of a regular $n$-gon in the plane. 
It can be presented as $$D_{n}=<r,s \; |\;  r^n=s^2=e,\; s^{-1}rs=r^{-1}>.$$
Let $\theta=\frac{2\pi}{n}$, for $n$ even, the irreducible representations are listed in the following table. If $n$ is odd, the representation $\phi_3$ and $\phi_4$ are omitted, and $k$ ranges from $1$ to $\frac{n-1}{2}$. 
\begin{table}[h]
\centering
\renewcommand{\tablename}{Table}
\renewcommand{\arraystretch}{1.2}  \doublerulesep 2.0pt
\begin{tabular}{|c|c|c|}
\hline
  & $r^j$ & $r^js$ \\
\hline
$\phi_1$ & 1 & 1\\
\hline
$\phi_2$ & 1 & -1\\
\hline
$\phi_3$ & $(-1)^{j}$ & $(-1)^{j}$\\
\hline
$\phi_4$ & $(-1)^{j}$ & $(-1)^{j+1}$\\
\hline
\makecell{$\rho_k$ \\ $k=1,\dots,\frac{n}{2}-1$}  & $ \begin{pmatrix}
 \cos{kj\theta} & -\sin{kj\theta} \\
 \sin{kj\theta} & \cos{kj\theta}
\end{pmatrix} $ & 
$ \begin{pmatrix}
 \cos{kj\theta} & \sin{kj\theta} \\
 \sin{kj\theta} & -\cos{kj\theta}
 \end{pmatrix} $ \\
\hline
\end{tabular}
\caption{The irreducible group representations for $D_n$}
\label{action}
\end{table}
\\
Given any representation $\mathscr{D}$ of $D_{n}$,
$$\mathscr{D}: D_n \rightarrow GL_k(\mathbb{R})$$
the character $\chi$ is a complex-valued function on $D_n$ given by $\chi(A)=Tr(\mathscr{D}(A)),$ for any $A\in G.$ Given two representations $\mathscr{D}_1$ and $\mathscr{D}_2$, define the inner product of characters $\chi_1$ and $\chi_2$ of $\mathscr{D}_1$ and $\mathscr{D}_2$, respectively,  as following
$$(\chi_1,\chi_2)=\frac{1}{|D_n|}\sum_{A\in D_{n}}\overline{\chi_1}(A) \chi_2(A).$$
A classic result in group representations \cite{2012Representation} states that the characters for distinct irreducible representations are orthogonal. Moreover, any representation $\mathscr{D}$ for $D_{n}$ is equivalent to the direct sum of irreducible representations, i.e., there is an invertible matrix $P$ such that
$$P^{-1}\mathscr{D}P=n_1 \mathscr{D}_1 \oplus n_2 \mathscr{D}_2 \oplus \dots \oplus n_h \mathscr{D}_h$$
where $\mathscr{D}_i$ are irreducible representations and
$$n_i = (\chi,\chi_i)\in\mathbb{Z},~ i = 1, \dots ,h.$$

Place $n$ equal masses at the vertices of the regular $n$-gon. The configuration of the $n$ particles can be represented by a vector $(q_1, q_2, \ldots, q_n) \in \mathbb{R}^{2n}$. For each element $A \in D_n$, it acts on the configuration by a permutation of the particles. For example, the action of $s$ is the reflection and the action of $r$ is the rotation by $\theta$. Corresponding to the position vector, the action of the symmetry group $D_n$ is represented by $2n\times 2n$ matrices,  
$$\mathscr{D}(e)=I_{2n},$$
$$\mathscr{D}(s)=
\begin{pmatrix}
F     &    0       & \dots &  0 & 0 \\
0     &    0       & \dots & 0  & F \\
0     &    0       & \dots & F  & 0\\
\vdots   & \vdots     & \begin{sideways}$\ddots$\end{sideways} & 0 & 0\\
0  & F &\dots & 0 & 0
\end{pmatrix}, \quad where \quad
F=\begin{pmatrix}
1 & 0 \\
0 & -1
\end{pmatrix},$$

$$\mathscr{D}(r)=
\begin{pmatrix}
0         &    0       & \dots &  0 & R(\theta) \\
R(\theta) &    0       & \dots & 0  & 0 \\
0         & R(\theta)  & \dots & 0  & 0\\
\vdots   & \vdots     & \ddots & 0 & 0\\
0  & 0 &\dots & R(\theta) & 0
\end{pmatrix}, \quad where \quad
R(\theta)=\begin{pmatrix}
\cos{\theta} & -\sin{\theta} \\
\sin{\theta} & \cos{\theta}
\end{pmatrix}.$$

The characters of $\mathscr{D}$ and the irreducible group representations can be easily evaluated. These numbers are listed in Table \ref{trace} (For $n$ odd, the representation $\phi_3$ and $\phi_4$ are omitted and $j$ goes from $1$ to $\frac{n-1}{2}$).

From the character table, our representation $\mathscr{D}$ can be decomposed into the direct sum of irreducible group representations. For $n$ even, 
$$(\chi_{\phi_i},\mathscr{D})=\frac{1}{2n}\times 2n=1,\quad(\chi_{\rho_i},\mathscr{D})=\frac{1}{2n}\times 2 \times 2n=2,$$
the representation $\mathscr{D}$ is equivalent to 
$$\phi_1 \oplus \phi_2 \oplus \phi_3 \oplus \phi_4 \oplus 2 \rho_1 \oplus \dots \oplus 2\rho_{\frac{n}{2}-1}.$$
For $n$ odd, the representation $\phi_3$ and $\phi_4$ are omitted and $\mathscr{D}$ is equivalent to
$$\phi_1 \oplus \phi_2  \oplus 2 \rho_1 \oplus \dots \oplus 2\rho_{\frac{n-1}{2}}.$$
\begin{table}[h]
\centering
\renewcommand{\tablename}{Table}
\renewcommand{\arraystretch}{1.2}  \doublerulesep 2.0pt
\begin{tabular}{|c|c|c|c|c|c|c|}
\hline
 A/ $\chi$ & $\phi_1$ & $\phi_2$ & $\phi_3$ & $\phi_4$ & $\rho_k$ & $\mathscr{D}$\\
\hline
$e$ & 1 & 1 & 1 & 1 & 2 & 2n \\
\hline
\makecell{$r^j,r^{n-j}$ \\ $j=1,2,\dots,\frac{n}{2}$} & 1 & 1 & $(-1)^j$ & $(-1)^j$ & $ 2\cos{kj\theta}$ & 0\\
\hline
$s,r^{2}s,\dots,r^{n-2}s$  & 1 & -1 & 1 & -1 & 0 & 0\\
\hline
$rs,r^{3}s,\dots,r^{n-1}s$ & 1 & -1 & -1 & 1 & 0 & 0\\
\hline
\end{tabular}
\caption{The trace of irreducible group representations for $D_n$}
\label{trace}
\end{table}

\subsection{The invariant subspace}
Let $H$ be a matrix which is invariant under $\mathscr{D}$, i.e., for every $a \in D_n$,
$$\mathscr{D}(a)H=H\mathscr{D}(a).$$
Suppose $E^{a}_{\lambda}$ is the eigenspace corresponding to $\lambda$ for the matrix $\mathscr{D}(a).$ For $v^a_{\lambda} \in E^{a}_{\lambda},$
$$\mathscr{D}(a)Hv^a_{\lambda}=H\mathscr{D}(a)v^a_{\lambda}=\lambda Hv^a_{\lambda}.$$
$$\Rightarrow Hv^a_{\lambda} \in E^{a}_{\lambda}$$
It shows the action $H$ does not mix the eigenspace for $\mathscr{D}(a)$, for any $a \in D_n$. For the same reason, the action $H$ does not mix the generalized eigenspaces of $\mathscr{D}(a)$ either. 

\

We want to find all the eigenvalues and eigenvectors of the Hessian $H$. or any matrix with $D_n$ symmetry. Since eigenvectors and eigenvalues of the irreducible representations of $\mathscr{D}(a)$ are easy to find, so instead, we can choose proper eigenvectors of $\mathscr{D}(a)$ as a new basis for $\mathbb{R}^{2n}$ and arrange them in such a way that all the matrices $\mathscr{D}(a)$ can be simultaneously simplified as a block-diagonal matrix. By the invariance principle, this new basis also put $H$ into block diagonal form.

The problem is when $\mathscr{D}(a)$ has eigenvalues with high multiplicity, then the invariant subspace, the eigenspace, can be quite big. However, if this eigenspace intersect with an eigenspace of another matrix $\mathscr{D}(b)$, where $b \in D_n$, then the invariant subspace can be further divided into smaller spaces. What we do next is to find the smallest invariant subspaces by analysing intersections of eigenspaces with different group elements.

We first discuss the cases with $n$ even. Since $\mathscr{D}$ is equivalent to 
$$\phi_1 \oplus \phi_2 \oplus \phi_3 \oplus \phi_4 \oplus 2 \rho_1 \oplus \dots \oplus 2\rho_{\frac{n}{2}-1},$$
there exists an invertible matrix $P$ such that
$$P^{-1}\mathscr{D}(s)P=
\begin{pmatrix}
1 & 0 & 0 & 0 & \dots & 0 & 0 \\
0 &-1 & 0 & 0 & \dots & 0 & 0\\
0 & 0 & 1 & 0 & \dots & 0 & 0\\
0 & 0 & 0 & -1 &\dots & 0 & 0\\
\vdots & \vdots & \vdots & \vdots & \ddots & \vdots & \vdots \\
0 & 0 & 0 & 0 & \dots & 1 & 0 \\
0 & 0 & 0 & 0 & \dots & 0 &-1 \\
\end{pmatrix}$$
$$P^{-1}\mathscr{D}(r)P=
\begin{pmatrix}
1 & 0 & 0 & 0 & 0 & 0 & \dots & 0 & 0 \\
0 & 1 & 0 & 0 & 0 & 0 & \dots & 0 & 0\\
0 & 0 & -1 & 0 & 0 & 0 & \dots & 0 & 0\\
0 & 0 & 0 & -1 & 0 & 0 & \dots & 0 & 0\\
0 & 0 & 0 & 0 & R(\theta) & 0 & \dots & 0 & 0 \\
0 & 0 & 0 & 0 & 0 & R(\theta)  & \dots & 0 & 0 \\
\vdots & \vdots & \vdots & \vdots & \vdots & \vdots & \ddots & \vdots & \vdots \\
0 & 0 & 0 & 0 & 0 & 0 & \dots & R((\frac{n}{2}-1)\theta) & 0 \\
0 & 0 & 0 & 0 & 0 & 0 & \dots & 0 & R((\frac{n}{2}-1)\theta) \\
\end{pmatrix}$$

\

The matrix $\mathscr{D}(s)$ has two eigenvalues: $1$ and $-1$, each with multiplicity $n$. 
Let $e_i, i=1,\dots,2n$ be the standard basis vectors in $R^{2n}$ with $(e_i)_j=\delta_{ij}.$
The eigenvalue $1$ has eigenvectors
$$v_1^{1,s}=e_1, \quad v_{2}^{1,s}=e_{n+1},$$
$$v_{2i-1}^{1,s}=e_{2i-1}+e_{2n+3-2i},\quad v_{2i}^{1,s}=e_{2i}-e_{2n+4-2i},$$
where $i$ ranges from $2$ to $\frac{n}{2}.$ And for eigenvalue $-1$, eigenvectors are 
$$v_1^{-1,s}=e_2,\quad v_{2}^{-1,s}=e_{n+2},$$
$$ v_{2i-1}^{-1,s}=e_{2i-1}-e_{2n+3-2i},\quad v_{2i}^{-1,s}=e_{2i}+e_{2n+4-2i},$$
where $i$ ranges from $2$ to $\frac{n}{2}.$ 

\

Eigenvalues for $\mathscr{D}(r)$ are $1,-1,e^{\pm ij\theta}, j=1,\dots, \frac{n}{2}-1.$ Eigenvalue $1$ has eigenvectors
$$v_1^{1,r}=\frac{1}{\sqrt{n}}\begin{pmatrix}
1 & 0 & \dots & \cos{k\theta} & \sin{k\theta} & \dots & \cos{(n-1)\theta} & \sin{(n-1)\theta}
\end{pmatrix}^{T},$$
$$v_2^{1,r}=\frac{1}{\sqrt{n}}\begin{pmatrix}
0 & 1 & \dots & -\sin{k\theta} & \cos{k\theta} & \dots & -\sin{(n-1)\theta} & \cos{(n-1)\theta}
\end{pmatrix}^{T}.$$

Eigenvalue $-1$ also has two eigenvectors
$$v_1^{-1,r}=\frac{1}{\sqrt{n}}\begin{pmatrix}
1 & 0 & \dots & (-1)^{k}\cos{k\theta} & (-1)^{k}\sin{k\theta} & \dots & -\cos{(n-1)\theta} & -\sin{(n-1)\theta}
\end{pmatrix}^{T},$$
$$v_2^{-1,r}=\frac{1}{\sqrt{n}}\begin{pmatrix}
0 & 1 & \dots & (-1)^{k+1}\sin{k\theta} & (-1)^{k}\cos{k\theta} & \dots & \sin{(n-1)\theta} & -\cos{(n-1)\theta}
\end{pmatrix}^{T}.$$

For eigenvalue $e^{ij\theta},$ its eigenvectors are
$$v_1^{j\theta,r}=\sqrt{\frac{2}{n}}\begin{pmatrix}
1 & 0 & \dots & e^{-ijk\theta}\cos{k\theta} & e^{-ijk\theta}\sin{k\theta} & \dots & e^{-ij(n-1)\theta}\cos{(n-1)\theta} & e^{-ij(n-1)\theta}\sin{(n-1)\theta}
\end{pmatrix}^{T},$$
$$v_2^{j\theta,r}=\sqrt{\frac{2}{n}}\begin{pmatrix}
0 & 1 & \dots & -e^{-ijk\theta}\sin{k\theta} & e^{-ijk\theta}\cos{k\theta} & \dots & -e^{-ij(n-1)\theta}\sin{(n-1)\theta} & e^{-ij(n-1)\theta}\cos{(n-1)\theta}
\end{pmatrix}^{T}.$$
\

Observe that 
$$E^s_1\cap E^r_1=Span\{v^{1,r}_1\},$$
$$E^s_{-1}\cap E^r_1=Span\{v^{1,r}_2\},$$
$$E^s_1\cap E^r_{-1}=Span\{v^{-1,r}_1\},$$
$$E^s_{-1}\cap E^r_{-1}=Span\{v^{-1,r}_2\},$$
$$E^s_{1}\cap \{E^r_{j\theta}+E^r_{-j\theta}\}=Span\{Re(v_1^{j\theta,r}),Im(v_2^{j\theta,r}) \},$$
$$E^s_{-1}\cap \{E^r_{j\theta}+E^r_{-j\theta}\}=Span\{Re(v_2^{j\theta,r}),Im(v_1^{j\theta,r}) \}.$$

By the invariance principle, the action $H$ does not mix these subspaces in the above. Choosing $v^{1r}_1, v^{1r}_2, ,v^{-1r}_1, v^{-1r}_2,Re(v_1^{j\theta,r}),$
$Im(v_2^{j\theta,r}),Re(v_2^{j\theta,r}),Im(v_1^{j\theta,r}), j=1,\dots,\frac{n}{2}-1$ as a new basis, the matrix $H$ is given by
\begin{align}\label{ma}
\begin{pmatrix}
\lambda_1 & 0 & 0 & 0 & 0 & 0 & \dots & 0 & 0 \\
0 & \lambda_2 & 0 & 0 & 0 & 0 & \dots & 0 & 0\\
0 & 0 & \lambda_3 & 0 & 0 & 0 & \dots & 0 & 0\\
0 & 0 & 0 & \lambda_4 & 0 & 0 & \dots & 0 & 0\\
0 & 0 & 0 & 0 & A_1 & 0 & \dots & 0 & 0 \\
0 & 0 & 0 & 0 & 0 & B_1  & \dots & 0 & 0 \\
\vdots & \vdots & \vdots & \vdots & \vdots & \vdots & \ddots & \vdots & \vdots \\
0 & 0 & 0 & 0 & 0 & 0 & \dots & A_{\frac{n}{2}-1} & 0 \\
0 & 0 & 0 & 0 & 0 & 0 & \dots & 0 & B_{\frac{n}{2}-1}
\end{pmatrix},    
\end{align}
where $A_i,B_i,i=1,\dots,\frac{n}{2}-1$ are $2\times 2$ matrix. 

For $n$ odd, the representation $\mathscr{D}$ is equivalent to
$$\phi_1 \oplus \phi_2  \oplus 2 \rho_1 \oplus \dots \oplus 2\rho_{\frac{n-1}{2}}.$$
There exists an invertible matrix $P$ such that
$$P^{-1}\mathscr{D}(s)P=
\begin{pmatrix}
1 & 0 & 0 & 0 & \dots & 0 & 0 \\
0 &-1 & 0 & 0 & \dots & 0 & 0\\
0 & 0 & 1 & 0 & \dots & 0 & 0\\
0 & 0 & 0 & -1 &\dots & 0 & 0\\
\vdots & \vdots & \vdots & \vdots & \ddots & \vdots & \vdots \\
0 & 0 & 0 & 0 & \dots & 1 & 0 \\
0 & 0 & 0 & 0 & \dots & 0 &-1 \\
\end{pmatrix}$$
$$P^{-1}\mathscr{D}(r)P=
\begin{pmatrix}
1 & 0 & 0 &  0 & \dots & 0 & 0 \\
0 & 1 & 0 & 0 & \dots & 0 & 0\\
0 & 0 & R(\theta) & 0 & \dots & 0 & 0 \\
0 & 0 & 0 & R(\theta)  & \dots & 0 & 0 \\
\vdots & \vdots & \vdots & \vdots & \ddots & \vdots & \vdots \\
0 & 0 & 0 & 0 & \dots & R((\frac{n-1}{2})\theta) & 0 \\
0 & 0 & 0 & 0 & \dots & 0 & R((\frac{n-1}{2})\theta) \\
\end{pmatrix}$$
Similarly, the matrix $\mathscr{D}(s)$ has two eigenvalues: $1$ and $-1$, each with multiplicity $n$, while eigenvalues for $\mathscr{D}(r)$ are $1,e^{\pm ij\theta}, j=1,\dots, \frac{n-1}{2}.$ Then it has the invariant space with
$$E^s_1\cap E^r_1=Span\{v^{1,r}_1\},$$
$$E^s_{-1}\cap E^r_1=Span\{v^{1,r}_2\},$$
$$E^s_{1}\cap \{E^r_{j\theta}+E^r_{-j\theta}\}=Span\{Re(v_1^{j\theta,r}),Im(v_2^{j\theta,r}) \},$$
$$E^s_{-1}\cap \{E^r_{j\theta}+E^r_{-j\theta}\}=Span\{Re(v_2^{j\theta,r}),Im(v_1^{j\theta,r}) \},$$
where $j$ ranges from $1$ to $\frac{n-1}{2}$. Choosing corresponding vectors as a new basis. It simplifies $H$ as a block-diagonal matrix as (\ref{ma}) with $\lambda_3,\lambda_4$ omitted and $i$ goes from $1$ to $\frac{n-1}{2}$ for $A_i,B_i.$ This greatly reduces the complexity of any matrix $H$ that's invariant under the group action, its eigenvalues can then be simply calculated.

\section{Application}

In this section, we apply the results of the previous section to several specific examples. In section 3.1, we consider the central configurations of the $n$-body problem with $n$ equal masses placed at a regular $n$-gon in the plane. We compute the eigenvalues of the Hessians. In section 3.2, we consider we consider the central configurations of the $1+n$-body problem with $n$ equal masses placed at a regular $n$-gon in the plane, and one with different mass at the center of the $n$-gon. In section 3.3, we consider the stability problem of the relative equilibria with the regular $n$-gon configuration studied in Section 3.1.

\subsection{The regular $n-$gon configuration}
 Consider the $n$ body problem with equal masses, let 
$$I=\frac{1}{2}\sum_{i=1}^{n}|q_i|^2,\quad U=\sum_{1 \le i<j \le n} \frac{1}{|q_i-q_j|}.$$
$I$ is the moment of inertia and $U$ is the potential function of the $n$-body system.  Let $q_k=(\cos{(k-1)}\theta,\sin{(k-1)\theta}),k=1,\dots,n.$ It is a regular $n$-gon configuration for the $n$ body problem which has $D_n$ as the symmetric group. The configuration $z_0=(q_1,\dots,q_n)$ is a critical point of $f=\sqrt{I}U.$ Let $z=(x_1,y_1,\dots,x_n,y_n)$ and $H=D^2f(z_0).$

Observe
$$U(\mathscr{D}(s)z)=U(x_1,-y_1,x_n,y_n,\dots,x_2,-y_2)=U(z).$$
Differential the equation with $z$
$$\mathscr{D}(s)'\nabla U(\mathscr{D}(s)z)=\nabla U(z),$$
where $\nabla U$ is the gradient of $U$ and $\mathscr{D}(s)'$ is the transpose of $\mathscr{D}(s)$.
Continue to differential with $z,$
$$\mathscr{D}(s)' D^2U(\mathscr{D}(s)z)\mathscr{D}(s)=D^2U(z).$$
As $\mathscr{D}(s)z_0=z_0$ and $\mathscr{D}(s)'=(\mathscr{D}(s))^{-1}=\mathscr{D}(s),$ it gives the invariant
$$D^2U(z_0)\mathscr{D}(s)=\mathscr{D}(s)D^2U(z_0).$$
Similarly, it can be verified
$$D^2U(z_0)\mathscr{D}(r)=\mathscr{D}(r)D^2U(z_0).$$
As the function $I$ also satisfies 
$$I(\mathscr{D}(s)z)=I(x_1,-y_1,x_n,y_n,\dots,x_2,-y_2)=I(z),$$
the Hessian $D^2(\sqrt{I})(z_0)$ is invariant under the action $\mathscr{D}$.
For
$$D^2f=\nabla(\sqrt{I})(\nabla U)'+UD^2(\sqrt{I})+\nabla U\nabla(\sqrt{I})'+\sqrt{I}D^2U,$$
the Hessian $H=D^2f(z_0)$ is invariant under the actions $\mathscr{D}$, i.e.,
$$H\mathscr{D}(s)=\mathscr{D}(s)H,\quad H\mathscr{D}(r)=\mathscr{D}(r)H.$$
As $\mathscr{D}(r^i)=\mathscr{D}(r)^i,$ we further have $\mathscr{D}(r^i)H=H\mathscr{D}(r^i).$
By the invariance, the following equations hold
\begin{equation}\label{is1}
  \frac{\partial f(z_0)}{\partial x_{n+2-k} \partial x_s}  = \frac{\partial f(z_0)}{\partial x_{k} \partial x_{n+2-s}},   \\
  -\frac{\partial f(z_0)}{\partial y_{n+2-k} \partial x_s} = \frac{\partial f(z_0)}{\partial y_{k} \partial x_{n+2-s}}, 
\end{equation}
\begin{equation}\label{is2}
\frac{\partial f(z_0)}{\partial x_{n+2-k} \partial y_s} = -\frac{\partial f(z_0)}{\partial x_{k} \partial y_{n+2-s}},  -\frac{\partial f(z_0)}{\partial y_{n+2-k} \partial y_s} =-\frac{\partial f(z_0)}{\partial y_{k} \partial y_{n+2-s}}.
\end{equation}
\begin{equation}\label{ir1}
    \cos{i\theta}\frac{\partial^2 f}{\partial x_{n+k-i} \partial x_s}-\sin{i\theta}\frac{\partial^2 f}{\partial y_{n+k-i} \partial x_s}=\cos{i\theta}\frac{\partial^2 f}{\partial x_k \partial x_{s+i}}+\sin{i\theta}\frac{\partial^2 f}{\partial x_k \partial y_{s+i}},
\end{equation}
\begin{equation}\label{ir2}
    \cos{i\theta}\frac{\partial^2 f}{\partial x_{n+k-i} \partial y_s}-\sin{i\theta}\frac{\partial^2 f}{\partial y_{n+k-i} \partial y_s}=-\sin{i\theta}\frac{\partial^2 f}{\partial x_k \partial x_{s+i}}+\cos{i\theta}\frac{\partial^2 f}{\partial x_k \partial y_{s+i}}, 
\end{equation}
\begin{equation}\label{ir3}
\sin{i\theta}\frac{\partial^2 f}{\partial x_{n+k-i} \partial x_s}+\cos{i\theta}\frac{\partial^2 f}{\partial y_{n+k-i} \partial x_s}=\cos{i\theta}\frac{\partial^2 f}{\partial y_k \partial x_{s+i}}+\sin{i\theta}\frac{\partial^2 f}{\partial y_k \partial y_{s+i}},
\end{equation} 
\begin{equation}\label{ir4}
\sin{i\theta}\frac{\partial^2 f}{\partial x_{n+k-i} \partial y_s}+\cos{i\theta}\frac{\partial^2 f}{\partial y_{n+k-i} \partial y_s}=-\sin{i\theta}\frac{\partial^2 f}{\partial y_k \partial x_{s+i}}+\cos{i\theta}\frac{\partial^2 f}{\partial y_k \partial y_{s+i}},
\end{equation}
For $n$ even, let $P=(v_1^T,\dots,v_{2n}^T)$ be a $2n \times 2n$ matrix defined by
\begin{small}
\begin{equation*}
(v^{1,r}_1,v^{1,r}_2,v^{-1,r}_1,v^{-1,r}_2,\dots,Re(v^{\theta,r}_{2k-1}),Im(v^{\theta,r}_{2k}),Im(v^{\theta,r}_{2k-1}),Re(v^{\theta,r}_{2k}),\dots),
\end{equation*}
\end{small}
where $k$ ranges from $1$ to $\frac{n}{2}-1$. For $n$ odd, we ignore $v^{-1,r}_1,v^{-1,r}_2$ and let
\begin{small}
\begin{equation*}
P=(v^{1,r}_1,v^{1,r}_2,\dots,Re(v^{\theta,r}_{2k-1}),Im(v^{\theta,r}_{2k}),Im(v^{\theta,r}_{2k-1}),Re(v^{\theta,r}_{2k}),\dots),
\end{equation*}
\end{small}
where $k$ ranges from $1$ to $\frac{n-1}{2}$. By the trigonometric identities
\begin{align} \label{tri}
\sum^{n}_{k=1}\sin{k\theta}=Im(\sum^{n}_{k=1}e^{ik\theta})=\frac{\sin{\frac{n\theta}{2}}\sin{\frac{(n+1)\theta}{2}}}{\sin{\frac{\theta}{2}}}, \\
\sum^{n}_{k=1}\cos{k\theta}=Re(\sum^{n}_{k=1}e^{ik\theta})=\frac{\sin{\frac{n\theta}{2}}\cos{\frac{(n+1)\theta}{2}}}{\sin{\frac{\theta}{2}}},    \nonumber
\end{align}
it is easy to show that the column vectors of $P$ are orthogonal and $P^{-1}=P^T.$ Then
$$P^TH=
\begin{pmatrix}
v_1\cdot(\frac{\partial f}{\partial z})_{x_1} & v_1\cdot(\frac{\partial f}{\partial z})_{y_1} & \dots &
v_1\cdot(\frac{\partial f}{\partial z})_{x_n} &
v_1\cdot(\frac{\partial f}{\partial z})_{y_n} \\
v_2\cdot(\frac{\partial f}{\partial z})_{x_1} & v_2\cdot(\frac{\partial f}{\partial z})_{y_1} & \dots &
v_2\cdot(\frac{\partial f}{\partial z})_{x_n} &
v_2\cdot(\frac{\partial f}{\partial z})_{y_n} \\
\vdots & \vdots & \ddots & \vdots & \vdots  \\
v_{2n}\cdot(\frac{\partial f}{\partial z})_{x_1} & v_{2n}\cdot(\frac{\partial f}{\partial z})_{y_1} & \dots &
v_{2n}\cdot(\frac{\partial f}{\partial z})_{x_n} &
v_{2n}\cdot(\frac{\partial f}{\partial z})_{y_n} \\
\end{pmatrix}$$
 Detailed computation (as shown in Appendix A), shows that $P^THP$ has the form as (\ref{ma}). And
 \begin{align*}
    \lambda_1=(P^THP)_{11}=&\frac{1}{\sqrt{n}}\sum_{i=1}^{n}\Big \{\cos{(i-1)\theta}\cdot v_1\cdot(\frac{\partial f}{\partial z})_{x_i}+\sin{(i-1)\theta}\cdot v_1\cdot(\frac{\partial f}{\partial z})_{y_i}\Big \} \\
    =&\frac{1}{n}\sum_{i=1}^{n}\Big \{\cos{(i-1)\theta}\sum_{j=1}^{n}(\cos{(j-1)\theta}\frac{\partial^2 f}{\partial x_j \partial x_i}+\sin{(j-1)\theta}\frac{\partial^2 f}{\partial y_j \partial x_i})\\\
    &+\sin{(i-1)\theta}\sum_{j=1}^{n}(\cos{(j-1)\theta}\frac{\partial^2 f}{\partial x_j \partial y_i}+\sin{(j-1)\theta}\frac{\partial^2 f}{\partial y_j \partial y_i})\Big \}
\end{align*}
by equations (\ref{ir1}) and (\ref{ir3}),
\begin{align*}
    \lambda_1
    =&\frac{1}{n}\sum_{i=1}^{n}\Big \{\cos{(i-1)\theta}\sum_{j=1}^{n}(\cos{(j-1)\theta}\frac{\partial^2 f}{\partial x_{n+i+1-j} \partial x_1}-\sin{(j-1)\theta}\frac{\partial^2 f}{\partial y_{n+i+1-j} \partial x_1})\\
    &+\sin{(i-1)\theta}\sum_{j=1}^{n}(\sin{(j-1)\theta}\frac{\partial^2 f}{\partial x_{n+i+1-j} \partial x_1}+\cos{(j-1)\theta}\frac{\partial^2 f}{\partial y_{n+i+1-j} \partial x_1})\Big \}\\
    =&\frac{1}{n}\sum_{i=1}^{n}\sum_{j=1}^{n}\Big \{\sin{(j-i)\theta}\frac{\partial^2 f}{\partial x_{n+i+1-j} \partial x_1}+\cos{(j-i)\theta}\frac{\partial^2 f}{\partial y_{n+i+1-j} \partial x_1}\Big \}. 
\end{align*}
From equations (\ref{is1}) and (\ref{is2}), we have
\begin{align*}
    \frac{\partial^2 f}{\partial x_{n+i+1-j}\partial x_1}=\frac{\partial^2 f}{\partial x_{n+2-(j-i+1)}\partial x_1}=\frac{\partial^2 f}{\partial x_{n+j-i+1}\partial x_1} \\
    -\frac{\partial^2 f}{\partial y_{n+i+1-j}\partial x_1}=-\frac{\partial^2 f}{\partial y_{n+2-(j-i+1)}\partial x_1}=\frac{\partial^2 f}{\partial y_{n+j-i+1}\partial x_1}.
\end{align*}
Therefore, it gives
\begin{align*}
    \lambda_1=\frac{2}{n}\sum^{n-1}_{p=1}(\frac{\partial^2 f}{\partial x_{n+1-p} \partial x_1}\cos{p\theta}-\frac{\partial^2 f}{\partial y_{n+1-p} \partial x_1}\sin{p\theta})(n-p)+\frac{\partial^2 f}{\partial x_{1} \partial x_1}.
\end{align*}
Similarly, the diagonal matrices (\ref{ma}) can be evaluated with
\begin{align*}
    \lambda_2=\frac{2}{n}\sum^{n-1}_{p=1}(\frac{\partial^2 f}{\partial x_{n+1-p} \partial y_1}\sin{p\theta}+\frac{\partial^2 f}{\partial y_{n+1-p} \partial y_1}\cos{p\theta})(n-p)+\frac{\partial^2 f}{\partial y_{1} \partial y_1},
\end{align*}
\begin{align*}
    \lambda_3=\frac{2}{n}\sum^{n-1}_{p=1}(-1)^p(\frac{\partial^2 f}{\partial x_{n+1-p} \partial x_1}\cos{p\theta}-\frac{\partial^2 f}{\partial y_{n+1-p} \partial x_1}\sin{p\theta})(n-p)+\frac{\partial^2 f}{\partial x_{1} \partial x_1},
\end{align*}
\begin{align*}
    \lambda_4=\frac{2}{n}\sum^{n-1}_{p=1}(-1)^p(\frac{\partial^2 f}{\partial x_{n+1-p} \partial y_1}\sin{p\theta}+\frac{\partial^2 f}{\partial y_{n+1-p} \partial y_1}\cos{p\theta})(n-p)+\frac{\partial^2 f}{\partial y_{1} \partial y_1},
\end{align*}
\begin{align*}
    (A_k)_{11}=\frac{2}{n}\sum^{n-1}_{p=1}(\frac{\partial^2 f}{\partial x_{n+1-p} \partial x_1}\cos{p\theta}-\frac{\partial^2 f}{\partial y_{n+1-p} \partial x_1}\sin{p\theta})(n-p)\cos{pk\theta}+\frac{\partial^2 f}{\partial x_{1} \partial x_1},
\end{align*}
\begin{align*}
    (A_k)_{22}=\frac{2}{n}\sum^{n-1}_{p=1}(\frac{\partial^2 f}{\partial x_{n+1-p} \partial y_1}\sin{p\theta}+\frac{\partial^2 f}{\partial y_{n+1-p} \partial y_1}\cos{p\theta})(n-p)\cos{pk\theta}+\frac{\partial^2 f}{\partial y_{1} \partial y_1},
\end{align*}
\begin{align*}
    (A_k)_{12}=(A_k)_{21}=\frac{2}{n}\sum^{n-1}_{p=1}(\frac{\partial^2 f}{\partial x_{n+1-p} \partial y_1}\cos{p\theta}-\frac{\partial^2 f}{\partial y_{n+1-p} \partial y_1}\sin{p\theta})(n-p)\sin{pk\theta},
\end{align*}
$$B_k=A_k.$$
Let $a=(A_k)_{11}+(A_k)_{22},b=|A_k|=(A_k)_{11}\cdot(A_k)_{22}-(A_k)_{12}^2.$ Then $H$ have eigenvalues with multiplicity 2
$$\lambda_{2k+3}=\frac{a+\sqrt{a^2-4b}}{2},\quad \lambda_{2k+4}=\frac{a-\sqrt{a^2-4b}}{2}, \quad k=1,\dots,\frac{n}{2}-1.$$
For $n$ odd, we omit eigenvalues $\lambda_3,\lambda_4$ and let $k$ goes from $1$ to $\frac{n-1}{2}$. Moreover, for the potential function $U$, as shown before, the Hessian $\frac{\partial^2U}{\partial z^2}$ at the configuration is also invariant under the dihedral group. Therefore the formulas for eigenvalues of the Hessian are still valid with $f=U$.

\subsection{The $n+1$ body problem}
For the $n+1$ body problem, let 
$$I=\frac{1}{2}\sum_{i=1}^{n}|q_i|^2+m|q_{n+1}|^2,\quad U=\sum_{1 \le i<j \le n} \frac{1}{|q_i-q_j|}+\sum_{i=1}^n \frac{m}{|q_i-q_{n+1}|}.$$
Consider the central configuration with $n$ equal masses placed at the vertices of a regular n-gon and mass $m$ placed at the center of gravity. Similarly let $q_k=(\cos{(k-1)}\theta,\sin{(k-1)\theta}),k=1,\dots,n$ and $q_{n+1}=(0,0).$ Then $z_0=(q_1,q_2,\dots,q_{n+1})$ is a central $+$ regular $n-$gon configuration. It also exhibits the dihedral symmetry with action $s$ as reflection and $r$ as rotation by $\theta$. Let the Hessian of the function $\sqrt{I}U$ at the configuration be $H.$ For group actions acting on the configuration in $R^{2n+2}$, define the representation with degree $2n+2$
$$\mathscr{D}(e)=I_{2n+2}$$
$$\mathscr{D}(s)=
\begin{pmatrix}
F     &    0       & \dots &  0 & 0 & 0 \\
0     &    0       & \dots & 0  & F & 0\\
0     &    0       & \dots & F  & 0 & 0\\
\vdots   & \vdots     & \begin{sideways}$\ddots$\end{sideways} & 0 & 0\\
0  & F &\dots & 0 & 0 & 0\\
0     &    0       & \dots &  0 & 0 & F
\end{pmatrix}, \quad where \quad
F=\begin{pmatrix}
1 & 0 \\
0 & -1
\end{pmatrix}$$

$$\mathscr{D}(r)=
\begin{pmatrix}
0         &    0       & \dots &  0 & R(\theta) & 0\\
R(\theta) &    0       & \dots & 0  & 0 & 0\\
0         & R(\theta)  & \dots & 0  & 0 & 0\\
\vdots   & \vdots     & \ddots & 0 & 0 & 0\\
0  & 0 &\dots & R(\theta) & 0 & 0 \\
0         &    0       & \dots &  0 & 0 & R(\theta)\\
\end{pmatrix}, \quad where \quad
R(\theta)=\begin{pmatrix}
\cos{\theta} & -\sin{\theta} \\
\sin{\theta} & \cos{\theta}
\end{pmatrix}.$$
For $n$ even, the representation is equivalent to 
$$\phi_1 \oplus \phi_2 \oplus \phi_3 \oplus \phi_4 \oplus 3 \rho_1 \oplus 2\rho_2 \oplus \dots \oplus 2\rho_{\frac{n}{2}-1}.$$
Similarly, it can be verified that $H$ is invariant under the group action, i.e.,
$$\mathscr{D}H=H\mathscr{D}.$$
Let $\Tilde{v}^{i,j}_k=(v^{i,j}_k,0,0)^T \in R^{2n+2}$ with $v^{i,j}_k$ are eigenvalues defined in Section $3.2$. They are eigenvectors for $\mathscr{D}$. It has a similar invariant space 
\begin{gather*}
E^s_1\cap E^r_1=Span\{\Tilde{v}^{1,r}_1\},\quad E^s_{-1}\cap E^r_1=Span\{\Tilde{v}^{1,r}_2\},\\
E^s_1\cap E^r_{-1}=Span\{\Tilde{v}^{-1,r}_1\}, \quad E^s_{-1}\cap E^r_{-1}=Span\{\Tilde{v}^{-1,r}_2\},\\
E^s_{1}\cap \{E^r_{\theta}+E^r_{-\theta}\}=Span\{Re(\Tilde{v}_1^{\theta,r}),Im(\Tilde{v}_2^{\theta,r}), e_{2n+1}\},\\
E^s_{-1}\cap \{E^r_{\theta}+E^r_{-\theta}\}=Span\{Re(\Tilde{v}_2^{\theta,r}),Im(\Tilde{v}_1^{\theta,r}), e_{2n+2}\},\\
    E^s_{1}\cap \{E^r_{j\theta}+E^r_{-j\theta}\}=Span\{Re(\Tilde{v}_1^{j\theta,r}),Im(\Tilde{v}_2^{j\theta,r})\},\\
    E^s_{-1}\cap \{E^r_{j\theta}+E^r_{-j\theta}\}=Span\{Re(\Tilde{v}_2^{j\theta,r}),Im(\Tilde{v}_1^{j\theta,r})\},
\end{gather*}
where $j$ ranges from $2$ to $\frac{n}{2}-1.$ By the invariance, choosing the corresponding vectors as a new basis, the matrix will be in a block-diagonal form as following
\begin{align}
\begin{pmatrix}
\lambda_1 & 0 & 0 & 0 & 0 & 0 & \dots & 0 & 0 \\
0 & \lambda_2 & 0 & 0 & 0 & 0 & \dots & 0 & 0 \\
0 & 0 & \lambda_3 & 0 & 0 & 0 & \dots & 0 & 0 \\
0 & 0 & 0 & \lambda_4 & 0 & 0 & \dots & 0 & 0 \\
0 & 0 & 0 & 0 & A_1 & 0 & \dots & 0 & 0\\
0 & 0 & 0 & 0 & 0 & B_1  & \dots & 0 & 0 \\
\vdots & \vdots & \vdots & \vdots & \vdots & \vdots & \ddots & \vdots & \vdots \\
0 & 0 & 0 & 0 & 0 & 0 & \dots & A_{\frac{n}{2}-1} & 0 \\
0 & 0 & 0 & 0 & 0 & 0 & \dots & 0 & B_{\frac{n}{2}-1} \\
\end{pmatrix},    
\end{align}
where $A_1,B_1$ are $3\times3$ matrix and $A_i,B_i, i=2,\dots,\frac{n}{2}-1$ are $2\times2$ matrix. 

\

For $n$ odd, the representation of $D_n$ is equivalent to 
$$\phi_1 \oplus \phi_2 \oplus 3 \rho_1 \oplus 2\rho_2 \oplus \dots \oplus 2\rho_{\frac{n-1}{2}}.$$
The smallest invariant subspaces are same as before, except that we need to drop the following
\begin{gather*}
    E^s_1\cap E^r_{-1}=Span\{\Tilde{v}^{-1,r}_1\}, \quad E^s_{-1}\cap E^r_{-1}=Span\{\Tilde{v}^{-1,r}_2\}.
\end{gather*}
Choosing corresponding vectors as a new basis, the Hessian $H$ will be a block-diagonal matrix with $\lambda_3,\lambda_4$ omitted and $i$ goes from $1$ to $\frac{n-1}{2}$.

\subsection{Dynamics near relative equilibria} For the Newtonian $n$-body problem, for each planar central configuration $z_0 \in \mathbb{R}^{2n}$, i.e., a critical point of $f= \sqrt{I}U$, there is correspondingly a class of {\em relative equilibria}\/, where $n$ particles move in concentric circles, keeping their relative positions. To study the dynamics near these relative equilibria, it is convenient to use a rotational coordinate frame. For a relative equilibrium with constant rotational frequency $\omega$, we make a canonical change of coordinates
$$q_i=exp(-\omega J_1t)\xi_i,$$
$$p_i=exp(-\omega J_1t)\eta_i,$$
where
$$J_1=\begin{pmatrix}
    0 & 1 \\
    -1 & 0 \\
\end{pmatrix}.
$$
and 
$$exp(-\omega J_1t)=\begin{pmatrix}
    \cos{\omega t} & -\sin{\omega t} \\
    \sin{\omega t} & \cos{\omega t}.
\end{pmatrix}$$
Under the uniform rotating coordinate, the relative equilibria becomes a rest point. The new Hamiltonian becomes
$$H=\sum_{i=1}^n(|\eta_i|^2-\omega \xi_i^TJ_1\eta_i)-\sum_{1\le i <j \le n }\frac{1}{|\xi_i-\xi_j|}.$$
The motion corresponding to equations (\ref{h1}) becomes
\begin{equation*}
\begin{aligned}
\dot{\xi_i} &=\omega J\xi_i +\eta_i, \\
\dot{\eta_i} &=\omega J\eta_i+\sum_{\substack{j=1 \\ j \ne i}
}^{n}\frac{\xi_j-\xi_i}{|\xi_j-\xi_i|^3}.
\end{aligned}
\end{equation*}
The corresponding Hamiltonian equations can be written as second-order equations by eliminating $\eta_i$ as following
\begin{equation} \label{lin0}
\ddot{\xi}_i=2\omega J_1\dot{\xi}_i+\omega^2
\xi_i+\sum_{\substack{j=1 \\ j \ne i}
}^{n}\frac{\xi_j-\xi_i}{|\xi_j-\xi_i|^3}.
\end{equation}
To study local behavior near the equilibrium solution, we need to linearize the equation. Local stability and dynamical structures are largely determined by the eigenvalues of this linearization. There are two terms in the linearization. The term corresponding to the potential $$\sum_{\substack{j=1 \\ j \ne i}
}^{n}\frac{\xi_j-\xi_i}{|\xi_j-\xi_i|^3}$$ is the Hessian $H$ of the potential function $U$ at the regular $n-$gon. The rest of the terms from the second order equation are already linear. Let $\xi=(\xi_1,\dots,\xi_n)$. The linearization for equations (\ref{lin0}) are  
\begin{equation} \label{lin1}
\ddot{\xi}=2\omega J\dot{\xi}+\omega^2
\xi+D^2U(z_0)\xi,
\end{equation}
where $J$ is the block-diagonal matrix:
$$\begin{pmatrix}
    J_1 &  &   \\
        & \ddots & \\
        &         &J_1 \\
\end{pmatrix}.$$
As $U$ is invariant under the action of group $D_n$, our results from Section 3.2 show that its Hessian will be a block-diagonal matrix with small blocks, of the sizes at most $2\times 2.$ With the new basis defined in Section 3.1, define $\xi=P\zeta$. The linearized equation (\ref{lin0}) becomes
\begin{equation} \label{lin2}
\ddot{\zeta}=2\omega P^{-1}JP\dot{\zeta}+\omega^2
\zeta+P^{-1}D^2U(z_0)P\zeta.
\end{equation}

The eigenvalues for the linearization satisfy the following equation
\begin{equation} \label{sta}
    |(\lambda^2-\omega^2)E_{2n}-2\omega\lambda P^{-1}JP -P^{-1}D^2U(z_0)P|=0.
\end{equation}
 We already know the term $P^{-1}D^2U(z_0)P$ is in a block-diagonal matrix with small blocks. The problem here is that we need to block diagonalize all the linear terms simultaneously. For this, we need to analyze how invariant subspaces of the Hessian of the potential function behave under the action of matrix $J$. 
Observe that
$$Jv_1^{1,r}=-v_2^{1,r}, \quad Jv_2^{1,r}=v_1^{1,r},$$
$$Jv_1^{-1,r}=-v_2^{-1,r}, \quad Jv_2^{-1,r}=v_1^{-1,r},$$
$$Jv_1^{j\theta,r}=-v_2^{j\theta,r}, \quad Jv_2^{j\theta,r}=v_1^{j\theta,r}.$$
Moreover,
$$\mathscr{D}(s)Jv_1^{1,r}=-\mathscr{D}(s)v_2^{1,r}=v_2^{1,r}=-Jv_1^{1,r},$$
$$\mathscr{D}(r)Jv_1^{1,r}=-\mathscr{D}(r)v_2^{1,r}=-v_2^{1,r}=Jv_1^{1,r},$$
We see that
$$Jv_1^{1,r} \in E^s_{-1}\cap E^r_1.$$
Similarly, the actions of $J$ have the following property $$Jv_2^{1,r}\in E^s_{1}\cap E^r_1,$$
$$Jv_1^{-1,r} \in E^s_{-1}\cap E^r_{-1}, \quad Jv_2^{-1,r}\in E^s_{1}\cap E^r_{-1},$$
$$Jv_1^{j\theta,r},Jv_2^{j\theta,r},Jv_1^{-j\theta,r},Jv_2^{-j\theta,r} \in E^r_{j\theta}\cup E^r_{-j\theta}.$$
Therefore, the matrix $J$ takes vectors in $E^s_{-1}\cap E^r_1$ into $ E^s_{1}\cap E^r_1$, $ E^s_{1}\cap E^r_1$ to $E^s_{-1}\cap E^r_1$,   $E^s_{1}\cap E^r_{-1}$ to $ E^s_{-1}\cap E^r_{-1}$, and $E^s_{-1}\cap E^r_{-1}$ to $E^s_{1}\cap E^r_{-1}$. Moreover, the space $E^r_{j\theta}\cup E^r_{-j\theta}$ is invariant under $J$.
With vectors defined in Section $3.2$ as a new basis, the matrix $P^{-1}JP$  takes the following block diagonal form:  for $n$ even
$$\begin{pmatrix}
    A_1   &   &     &        & \\
        & A_2 &     &        & \\
        &   & B_1 &        &  \\
        &   &     & \ddots &   \\
        &   &     &        &  B_{\frac{n}{2}-1} 
\end{pmatrix},$$
where $A_k,k=1,2$ is an $2\times 2$ anti-diagonal matrix and $B_k, k=1,\dots,\frac{n}{2}-1$ is a $4\times4$ matrix.

When $n$ is odd, matrix $A_2$ is absent and $k$ ranges from $1$ to $\frac{n-1}{2}.$  

\

We conclude that the linearized equations at the regular $n-$gon relative equilibria can be simplified into block diagonal matrix with block sizes at most $4\times 4$. Therefore, the eigenvalue equation (\ref{sta}) with relative ease, comparing to the original equation.

\bibliographystyle{acm}
\bibliography{ref}
\newpage
\appendix
\section{}
We put some of the computations in Section 3.2 in this appendix.

\begin{align*}
    (P^THP)_{12}=&\frac{1}{\sqrt{n}}\sum_{i=1}^{n}-\sin{(i-1)\theta}\cdot v_1\cdot(\frac{\partial f}{\partial z})_{x_i}+\cos{(i-1)\theta}\cdot v_1\cdot(\frac{\partial f}{\partial z})_{y_i} \\
    =&\frac{1}{n}\sum_{i=1}^{n}-\sin{(i-1)\theta}\sum_{j=1}^{n}(\cos{(j-1)\theta}\frac{\partial^2 f}{\partial x_j \partial x_i}+\sin{(j-1)\theta}\frac{\partial^2 f}{\partial y_j \partial x_i})\\
    &+\cos{(i-1)\theta}\sum_{j=1}^{n}(\cos{(j-1)\theta}\frac{\partial^2 f}{\partial x_j \partial y_i}+\sin{(j-1)\theta}\frac{\partial^2 f}{\partial y_j \partial y_i})
\end{align*}
by equations (\ref{ir1}) and (\ref{ir3}),
\begin{align*}
    (P^THP)_{12}
    =&\frac{1}{n}\sum_{i=1}^{n}-\sin{(i-1)\theta}\sum_{j=1}^{n}(\cos{(j-1)\theta}\frac{\partial^2 f}{\partial x_{n+i+1-j} \partial x_1}-\sin{(j-1)\theta}\frac{\partial^2 f}{\partial y_{n+i+1-j} \partial x_1})\\
    &+\cos{(i-1)\theta}\sum_{j=1}^{n}(\sin{(j-1)\theta}\frac{\partial^2 f}{\partial x_{n+i+1-j} \partial x_1}+\cos{(j-1)\theta}\frac{\partial^2 f}{\partial y_{n+i+1-j} \partial x_1})\\
    =&\frac{1}{n}\sum_{i=1}^{n}\sum_{j=1}^{n}\sin{(j-i)\theta}\frac{\partial^2 f}{\partial x_{n+i+1-j} \partial x_1}+\cos{(j-i)\theta}\frac{\partial^2 f}{\partial y_{n+i+1-j} \partial x_1}. 
\end{align*}
From equations (\ref{is1}) and (\ref{is2}), we have
\begin{align*}
    (P^THP)_{12}
    =&\frac{1}{n}\sum_{i=1}^{n}\sum_{j=1}^{n}\sin{(j-i)\theta}\frac{\partial^2 f}{\partial x_{n+i+1-j} \partial x_1}+\cos{(j-i)\theta}\frac{\partial^2 f}{\partial y_{n+i+1-j} \partial x_1} \\
    =&\frac{1}{n}\sum_{j=1}^{n}\sum_{i=1}^{n}-\sin{(i-j)\theta}\frac{\partial^2 f}{\partial x_{n+j+1-i} \partial x_1}-\cos{(i-j)\theta}\frac{\partial^2 f}{\partial y_{n+j+1-i} \partial x_1} =- (P^THP)_{12},
\end{align*}
it gives $(P^THP)_{12}=0.$ Similarly, 
\begin{align*}
    (P^THP)_{13}
    =&\frac{1}{n}\sum_{i=1}^{n}\sum_{j=1}^{n}(-1)^{i-1}(\cos{(j-i)\theta}\frac{\partial^2 f}{\partial x_{n+i+1-j} \partial x_1}-\sin{(j-i)\theta}\frac{\partial^2 f}{\partial y_{n+i+1-j} \partial x_1}) \\
    =&\frac{1}{n}\sum_{p=0}^{n-1}\sum_{|j-i|=p}(-1)^{i-1}(\cos{p\theta}\frac{\partial^2 f}{\partial x_{n+1-p} \partial x_1}-\sin{p\theta}\frac{\partial^2 f}{\partial y_{n+1-p} \partial x_1}),
\end{align*}
and by equations (\ref{is1}) and (\ref{is2})
\begin{gather*}
  \frac{\partial f(z_0)}{\partial x_{n+1-p} \partial x_1}  = \frac{\partial f(z_0)}{\partial x_{n+2-(p+1)} \partial x_{1}}= \frac{\partial f(z_0)}{\partial x_{p+1} \partial x_{1}},   \\
  -\frac{\partial f(z_0)}{\partial y_{n+1-p} \partial x_1} = -\frac{\partial f(z_0)}{\partial y_{n+2-(p+1)} \partial x_{1}}=\frac{\partial f(z_0)}{\partial y_{p+1} \partial x_{1}}, \\
  \frac{\partial f(z_0)}{\partial y_{n+1-p} \partial y_1}= \frac{\partial f(z_0)}{\partial y_{n+2-(p+1)} \partial y_1}=\frac{\partial f(z_0)}{\partial y_{p+1} \partial y_{1}}. 
\end{gather*}
For $j-i=p,$ we have $i=1,2,\dots,n-p.$ For $i-j=p,$ we have $i=p+1,p+2,\dots,n.$ As $n$ is even, the sum $\sum_{|j-i|=p}(-1)^{i-1}$ is $\frac{(1+(-1)^{p})(1-(-1)^p)}{2}=0.$ Then we verify $(P^THP)_{13}=0.$
For
\begin{align*}
    (P^THP)_{14}
    =&\frac{1}{n}\sum_{i=1}^{n}\sum_{j=1}^{n}(-1)^{i-1}(\sin{(j-i)\theta}\frac{\partial^2 f}{\partial x_{n+i+1-j} \partial x_1}+\cos{(j-i)\theta}\frac{\partial^2 f}{\partial y_{n+i+1-j} \partial x_1}) \\
    =&\frac{1}{n}\sum_{p=0}^{n-1}(\sin{p\theta}\frac{\partial^2 f}{\partial x_{n+1-p} \partial x_1}+\cos{p\theta}\frac{\partial^2 f}{\partial y_{n+1-p} \partial x_1})(\frac{(1-(-1)^{p})(1-(-1)^p)}{2})\\
    =&\frac{1}{n}\sum_{p=0}^{n-1}a_p,
\end{align*}
Since $\frac{\partial^2 f}{\partial x_{n+1-p} \partial x_1}=\frac{\partial^2 f}{\partial x_{1+p} \partial x_1}, \frac{\partial^2 f}{\partial y_{n+1-p} \partial x_1}=-\frac{\partial^2 f}{\partial y_{1+p} \partial x_1}$, it gives $a_p=-a_{n-p}.$ As $a_0=0,$ we conclude $(P^THP)_{14}=0.$
For
\begin{align*}
    (P^THP)_{15}
    =&\frac{\sqrt{2}}{n}\sum_{i=1}^{n}\sum_{j=1}^{n}(\cos{(j-i)\theta}\frac{\partial^2 f}{\partial x_{n+i+1-j} \partial x_1}-\sin{(j-i)\theta}\frac{\partial^2 f}{\partial y_{n+i+1-j} \partial x_1})\cos{(i-1)\theta} \\
    =&\frac{\sqrt{2}}{n}\sum_{p=0}^{n-1}(\cos{p\theta}\frac{\partial^2 f}{\partial x_{n+1-p} \partial x_1}-\sin{p\theta}\frac{\partial^2 f}{\partial y_{n+1-p} \partial x_1})(\sum_{i=1}^{n-p}+\sum_{i=p+1}^{n})\cos{(i-1)\theta}\\
    =&\frac{\sqrt{2}}{n}\sum_{p=0}^{n-1}a_p,
\end{align*}
by the trigonometric identities (\ref{tri}),
\begin{align*}
  (\sum_{i=1}^{n-p}+\sum_{i=p+1}^{n})\cos{(i-1)\theta} &=1+(\sum_{i=1}^{n-p-1}+\sum_{i=1}^{n-1}-\sum_{i=1}^{p-1})\cos{i\theta} \\
  &=1+\frac{\sin{\frac{n-p-1}{2}\theta}\cos{\frac{n-p}{2}\theta}+\sin{\frac{n-1}{2}\theta}\cos{\frac{n}{2}\theta}-\sin{\frac{p-1}{2}\theta}\cos{\frac{p}{2}\theta}}{\sin{\frac{\theta}{2}}} \\
  &=-2\frac{\sin{\frac{p}{2}\theta}\cos{\frac{\theta}{2}}\cos{\frac{p}{2}\theta}}{\sin{\frac{\theta}{2}}} \\
  &=-\frac{\sin{p\theta}}{\tan{\frac{\theta}{2}}}.
\end{align*}
By $a_p=a_{n-p},a_0=0,$ we have $(P^THP)_{15}$ is zero. For
\begin{align*}
    (P^THP)_{35}
    =&\frac{\sqrt{2}}{n}\sum_{i=1}^{n}(v_3\cdot(\frac{\partial f}{\partial z})_{x_{i}}\cdot\cos{(i-1)\theta}+v_3\cdot(\frac{\partial f}{\partial z})_{y_{i}}\cdot\sin{(i-1)\theta})\cos{(i-1)\theta} \\
    =&\frac{\sqrt{2}}{n}\sum_{i=1}^{n}\sum_{j=1}^{n}(-1)^{j-1}\cos{(i-1)\theta}((\cos{(j-1)\theta}\frac{\partial^2 f}{\partial x_j \partial x_i}+\sin{(j-1)\theta}\frac{\partial^2 f}{\partial y_j \partial x_i})\cdot\cos{(i-1)\theta} \\
    &+(\cos{(j-1)\theta}\frac{\partial^2 f}{\partial x_j \partial y_i}+\sin{(j-1)\theta}\frac{\partial^2 f}{\partial y_j \partial y_i})\cdot\sin{(i-1)\theta})
\end{align*}
by equations (\ref{ir1}) and (\ref{ir3}),
\begin{align*}
    (P^THP)_{35}
     =&\frac{\sqrt{2}}{n}\sum_{i=1}^{n}\sum_{j=1}^{n}(-1)^{j-1}((\cos{(j-1)\theta}\frac{\partial^2 f}{\partial x_{n+i+1-j} \partial x_1}-\sin{(j-1)\theta}\frac{\partial^2 f}{\partial y_{n+i+1-j} \partial x_1})\cdot\cos{(i-1)\theta} \\
    &+(\sin{(j-1)\theta}\frac{\partial^2 f}{\partial x_{n+i+1-j} \partial x_1}+\cos{(j-1)\theta}\frac{\partial^2 f}{\partial y_{n+i+1-j} \partial x_1})\cdot\sin{(i-1)\theta})\cos{(i-1)\theta} \\
    =&\frac{\sqrt{2}}{n}\sum_{i=1}^{n}\sum_{j=1}^{n}(-1)^{j-1}((\cos{(j-i)\theta}\frac{\partial^2 f}{\partial x_{n+i+1-j} \partial x_1}-\sin{(j-i)\theta}\frac{\partial^2 f}{\partial y_{n+i+1-j} \partial x_1}))\cos{(i-1)\theta} \\
    =&\frac{\sqrt{2}}{n}\sum_{p=0}^{n-1}(\cos{p\theta}\frac{\partial^2 f}{\partial x_{n+1-p} \partial x_1}-\sin{p\theta}\frac{\partial^2 f}{\partial y_{n+1-p} \partial x_1})(\sum_{i=1}^{n-p}+\sum_{i=p+1}^{n})(-1)^{i+p-1}\cos{(i-1)\theta}\\
    =&\frac{\sqrt{2}}{n}\sum_{p=0}^{n-1}a_p.
\end{align*}
For
\begin{align*}
  (\sum_{i=1}^{n-p}+\sum_{i=p+1}^{n})(-1)^{i+p-1}\cos{(i-1)\theta}=&(-1)^p(\sum_{i=1}^{n-p}+\sum_{i=1}^{n}-\sum_{i=1}^{p})(-1)^{i-1}\cos{(i-1)\theta} \\
  =&(-1)^p(1+(\sum_{i=1}^{n-p-1}+\sum_{i=1}^{n-1}-\sum_{i=1}^{p-1})(-1)^i\cos{i\theta}) \\
  =&(-1)^p(1+(\sum_{i=1}^{n-p-1}+\sum_{i=1}^{n-1}-\sum_{i=1}^{p-1})\cos{i(\theta+\pi)}),
\end{align*}
by the trigonometric identities (\ref{tri}), we have 
\begin{align*}
  (\sum_{i=1}^{n-p}+\sum_{i=p+1}^{n})(-1)^{i+p-1}\cos{(i-1)\theta}=\sin{p\theta}\tan{\frac{\theta}{2}}.
\end{align*}
As $a_p=a_{n-p}$ and $a_0=0,$ $(P^THP)_{35}=0.$ For
\begin{align*}
    (P^THP)_{57}
    =&\sqrt{\frac{2}{n}}\sum_{i=1}^{n}(v_5\cdot(\frac{\partial f}{\partial z})_{x_{i}}\cdot\cos{(i-1)\theta}+v_3\cdot(\frac{\partial f}{\partial z})_{y_{i}}\cdot\sin{(i-1)\theta})\sin{(i-1)\theta} \\
    =&\frac{2}{n}\sum_{i=1}^{n}\sum_{j=1}^{n}[(\cos{(j-1)\theta}\frac{\partial^2 f}{\partial x_j \partial x_i}+\sin{(j-1)\theta}\frac{\partial^2 f}{\partial y_j \partial x_i})\cdot\cos{(i-1)\theta} \\
    &+(\cos{(j-1)\theta}\frac{\partial^2 f}{\partial x_j \partial y_i}+\sin{(j-1)\theta}\frac{\partial^2 f}{\partial y_j \partial y_i})\cdot\sin{(i-1)\theta}]\cos{(j-1)\theta}\sin{(i-1)\theta}
\end{align*}
by equations (\ref{ir1}) and (\ref{ir3}),
\begin{align*}
    (P^THP)_{57}
    =&\frac{\sqrt{2}}{n}\sum_{i=1}^{n}\sum_{j=1}^{n}((\cos{(j-i)\theta}\frac{\partial^2 f}{\partial x_{n+i+1-j} \partial x_1}-\sin{(j-i)\theta}\frac{\partial^2 f}{\partial y_{n+i+1-j} \partial x_1}))\cos{(j-1)\theta}\sin{(i-1)\theta} \\
    =&\frac{\sqrt{2}}{n}\sum_{p=0}^{n-1}(\cos{p\theta}\frac{\partial^2 f}{\partial x_{n+1-p} \partial x_1}-\sin{p\theta}\frac{\partial^2 f}{\partial y_{n+1-p} \partial x_1})b_p
\end{align*}
where $b_p=\sum_{j=1}^{n-p}\cos{(j-1)\theta}\sin{(j+p-1)\theta}+\sum_{j=p+1}^{n}\cos{(j-1)\theta}\sin{(j-p-1)\theta}$. As
\begin{align*}
    b_p=&\sum_{j=1}^{n-p}\cos{(j-1)\theta}(\sin{(j-1)\theta}\cos{p\theta}+\cos{(j-1)\theta}\sin{p\theta}) \\
    &+\sum_{j=p+1}^{n}\cos{(j-1)\theta}(\sin{(j-1)\theta}\cos{p\theta}-\cos{(j-1)\theta}\sin{p\theta}) \\
    =&(\sum_{j=1}^{n-p-1}+\sum_{j=1}^{n-1}-\sum_{j=1}^{p-1})\cos{j\theta}\sin{j\theta}\cos{p\theta}+(\sum_{j=1}^{n-p-1}-\sum_{j=1}^{n-1}+\sum_{j=1}^{p-1})\cos{j\theta}\cos{j\theta}\sin{p\theta}+\sin{p\theta}\\
    =&\frac{1}{2}[(\sum_{j=1}^{n-p-1}+\sum_{j=1}^{n-1}-\sum_{j=1}^{p-1})\sin{2j\theta}\cos{p\theta}+(\sum_{j=1}^{n-p-1}-\sum_{j=1}^{n-1}+\sum_{j=1}^{p-1})\cos{2j\theta}\sin{p\theta}+\sin{p\theta}]\\
    =&\frac{1}{2}(\cos{p\theta}\sin{2p\theta}+\sin{p\theta}\sin^2{p\theta})
\end{align*}
As $b_p=b_{n-p}$ and $b_0=0,$ $(P^THP)_{57}=0.$  By similar computation, it can be verified $P^THP$ has the form as (\ref{ma}).
\end{document}